\documentclass[12pt,twoside]{article}
\usepackage{amstex}
\usepackage{amssymb}
\usepackage{graphics}
\usepackage{epsf}
\newcommand{\proof}{\medskip \noindent{\bf Proof:\ \ }}

\newcommand{\remark}{\medskip \noindent{\bf Remark:\ \ }}
\begin{document}
\title{\bf{Weights of Markov Traces on Hecke Algebras}\footnote{1991 Subject Classification: 05E99}}
\author{R.C. Orellana}
\maketitle
\begin{abstract}
We compute the weights, i.e. the values at the minimal idempotents, for the Markov trace on the Hecke algebra of type $B_n$ and type $D_n$.  In order to prove the weight formula, we define representations of the Hecke algebra of type $B$ onto a reduced Hecke Algebra of type $A$.  To compute the weights for type $D$ we use the inclusion of the Hecke algebra of type $D$ into the Hecke algebra of type $B$. 
\end{abstract}
\pagestyle{myheadings}
\markboth{R. C. Orellana}{Weights of Markov Traces on Hecke Algebras}
\newtheorem{theorem}{Theorem}[section]
\newtheorem{definition}[theorem]{Definition}
\newtheorem{lemma}[theorem]{Lemma}
\newtheorem{claim}[theorem]{Claim}
\newtheorem{proposition}[theorem]{Proposition}
\newtheorem{cor}[theorem]{Corollary}
\centerline{\bf{Introduction}}
\vskip .07in
The Hecke algebra of type $B_n$, $H_n(q,Q)$, is semisimple whenever  $Q\neq -q^k$, $k\in \{0,\pm 1,\ldots \pm (n-1)\}$, and $q$ is  not a root of unity.  The simple components are indexed by pairs of Young diagrams.  Alternatively, these Hecke algebras can be defined as a finite dimensional quotient of the group algebra of the braid group of type $B$.  Motivated by their study of link invariants related to the braid group of type $B$, Geck and Lambropolou \cite{gl} have defined certain linear traces on the Hecke algebra of type $B$ called Markov traces.  Their definition is given inductively.
\vskip 0in
In this paper we give an alternative way of computing this trace.  Since the Hecke algebra of type $B$ is semisimple, any linear trace can be written as a weighted linear combination of the irreducible characters (the usual trace).  The coefficients in this linear expression are called weights.  The weights are equal to the values of the trace at the minimal idempotents. Since the characters are known, it follows that the weights completely determine the trace.  The weights are also indexed by pairs of Young diagrams.    
\vskip 0in
We have found the weight formula for the Markov trace defined by Geck and Lambropolou \cite{gl} for the Hecke algebra of type $B$. The weight formula can be written as a product of Schur functions. To prove this formula we construct a homomorphism from the specialization of the Hecke algebra of type $B$, $H_n(q,-q^{r_1+m})$, onto a reduced Hecke algebra of type $A$. Using this homomorphism we obtain that the Markov trace of the Hecke algebra of type $B$ appears as a pullback of the Markov trace of the reduced Hecke algebra of type $A$.       
\vskip 0in
A consequence of the above mentioned homomorphism is the existence of a duality between the quantum group $U_q(sl(r))$ and the specialized Hecke algebra $H_n(q,Q)$.  
\vskip 0in
In order to find the weights for the Hecke algebra of type $D$ we use the results of Hoefsmit \cite{h} on the inclusion of the Hecke algebra of type $D$ into the Hecke algebra of type $B$. We also use the results of Geck \cite{g} on obtaining Markov traces of the Hecke algebra of type $D$ from those of type $B$.     
\vskip 0in
In a future paper we will use the weight formula to determine the structure of a certain quotient of the Hecke algebra at roots of unity.  This quotient of the Hecke algebras is important in the theory of von Neumann algebras since they provide examples of subfactors of the  hyperfinite $II_1$ factor.
\vskip 0.1in
This paper is organized as follows.
In the first section we give some background definitions for partitions and Schur functions. We also define the braid groups of type $A$ and type $B$;  and we give representations of the braid group of type $B$ into the braid group of type $A$. In the second section we define the Hecke Algebra of type $B$. We also define the right coset representatives and distinguished double coset representatives of $H_n(q,Q)$ in $H_{n-1}(q,Q)$ which are needed in the later sections. \vskip 0in
In the third section we show that the representations given in section 1 for braid groups can be extended to well-defined representations of the Hecke algebra of type $B$ onto a reduced algebra of type $A$. 
In the fourth section we define the Markov trace and quote some of the results from \cite{gl} which are relevant in our quest. In the fifth section we prove the weight formula.  In the last section derive the weights for the Hecke algebra of type $D$ from the weights for the Hecke algebra of type $B$.  
\vskip 0.1in
\noindent
\emph{Acknowledgments: } This paper is a part of my Ph.D. thesis at University of California, San Diego.  I would like to thank my advisor Hans Wenzl for all the encouragement, for teaching me so much mathematics and for all the useful suggestions.
\section{Preliminaries and Notation}
\vskip 0in
Let $F$ be a field of characteristic $0$.  The algebra of $n\times n$ matrices over $F$ is denoted by  $M_n(F)=M_n$.  Let $A\subset B$ be an inclusion of semisimple algebras.  Then $A=\bigoplus_{i=1}^m A_i$ and $B=\bigoplus_{j=1}^lB_j$, with $A_i\cong M_{a_i}$, $B_j\cong M_{b_j}$ and $a_i,\ b_j\in \mathbb{N}$.  
\vskip 0in
Since a simple $B_j$ module  is also an $A$ module, it can be decomposed into a direct sum of simple $A$ modules. Let $g_{ij}$ be the number of simple $A_i$ modules in this decomposition.  The matrix $G=(g_{ij})$ is called the inclusion matrix for $A\subset B$. This matrix can be conveniently described by an inclusion diagram.  This is a graph with vertices arranged in 2 lines.   In one line, the vertices are in $1-1$ correspondence with the simple summands $A_i$ of $A$, in the other one with summands $B_j$ of $B$.  Then a vertex corresponding to $A_i$ is joined with a vertex corresponding to $B_j$ by $g_{ij}$ edges.  This is known as the Bratteli diagram.  
\vskip 0in
A trace is a linear funtional $tr:B\rightarrow F$ such that $tr(ab)=tr(ba)$ for all $a,b\in B$.  Recall that $tr$ is nondegenerate if for any $b\in B$, $b\neq 0$, there is some $b'\in B$ such that $tr(bb')\neq 0$.  There is up to scalar multiples only one trace on $M_n(F)$.  Any trace $tr$ on $B=\bigoplus B_j$ is completely determined by a vector $\vec{t}=(t_j)$, where $t_j=tr(p_j)$ and $p_j$ is a minimal idempotent of $B_j$.  The vector $\vec{t}$ is called the weight vector of $tr$.  In a similar manner we can characterize a trace on $A$ by a weight vector $\vec{s}=(s_i)$.  
The weight vectors $\vec{t}$ and $\vec{s}$ are related as follows:
$$G\vec{t}=\vec{s}.$$
\subsection{Partitions and Standard Tableaux}
\vskip 0.05in
A \emph{partition} of $n$ is a sequence of positive numbers $\alpha=(\alpha_1,\ldots, \alpha_k)$ such that $\alpha_1\geq \alpha_2\geq \ldots \geq \alpha_k \geq 0 $ and $n=|\alpha|=\alpha_1+\alpha_2+\ldots +\alpha_k$. The $\alpha_i$'s are called the \emph{parts} of $\alpha$, and the number of nonzero parts of $\alpha$ is called the \emph{length} of $\alpha$, denoted by $l(\alpha)$. If we say $l(\alpha)\leq r_1$, then we we will mean that there are $s\leq r_1$ nonzero parts and the remaining $r_1-s$ are equal to zero. We use the notation $\alpha\vdash n$ to mean $\alpha$ is a partition of $n$. 
\vskip 0in
A \emph{Young diagram} is a pictorial representation of a partition $\alpha$ as an array of $n$ boxes with $\alpha_1$ boxes in the first row, $\alpha_2$ boxes in the second row, and so on. We count the rows from top to bottom.  We shall denote the Young diagram and the partitions by the same symbol $\alpha$.  We use the word \emph{shapes} interchangeably with the word partitions.  The set of all partitions of $n$ is denoted by $\Lambda_n$.
\vskip 0in
If $\alpha$ and $\beta$ are two partitions with $|\alpha| \leq |\beta|$, then we write $\alpha\subset \beta$ if $\alpha_i\leq \beta_i$ for all $i$.  In this case we say that $\alpha$ is contained in $\beta$.  If $\alpha\subset \beta$ then the set-theoretic difference $\beta-\alpha$ is called a \emph{skew diagram}, $\alpha/\beta$.
\vskip 0in 
A \emph{standard tableaux} of shape $\alpha$ is a filling of the boxes with numbers $1, 2,\ldots, n$ such that the numbers in each row increase from left to right and in every column from top to bottom.  The notation $t^{\alpha}$ will be used to denote a standard tableaux of shape $\alpha$.   We say $t^{\alpha}\subset t^{\beta}$ if $t^{\alpha}$ is obtained by removing appropriate boxes from $t^{\beta}$, i.e. the numbers $1, 2, \ldots, |\alpha|$ are in the same boxes of both in $t^{\alpha}$ and in $t^{\beta}$. In this case we say that $t^{\alpha}$ is contained in $t^{\beta}$.
\vskip 0in
A \emph{double partition} of size $n$, $(\alpha,\beta)$, is an ordered pair of partitions $\alpha$ and $\beta$ such that $|\alpha|+|\beta|=n$.  The length of a double partition is $l(\alpha, \beta)=l(\alpha)+l(\beta)$.  A double partition can be associated with a pair of Young diagrams. The set of all pairs of Young diagrams with $n$ boxes is denoted by $\Lambda_n^D$.  We write as for single partitions $(\alpha,\beta) \subset (\delta,\gamma)$ if $\alpha_i \leq \delta_i$ and $\beta_j \leq \gamma_j$ for all $i$ and $j$.  That is, $(\alpha,\beta)$ is obtained from $(\delta,\gamma)$ by removing  appropriate boxes of $(\delta,\gamma)$, or equivalently $(\alpha,\beta)$ is contained in $(\delta,\gamma)$.   
\vskip 0in
$t^{(\alpha,\beta)}=(t^{\alpha},t^{\beta})$ is a pair of standard tableaux if the arrangement of the numbers is in increasing order in the rows and columns of both $t^{\alpha}$ and $t^{\beta}$.  For any double partition $(\alpha,\beta)$, let $T_{(\alpha,\beta)}$ denote the set of standard tableaux of shape $(\alpha,\beta)$ and $\mathcal{T}_n$ denote the set of all pairs of tableaux with $n$ boxes.  For example, the diagram below is a pair of standard tableaux of shape $([2,1,1],[3,2])$ 
\vskip 0in
\begin{picture}(200,60)(0,0)
\put(175,50){\framebox(12,12){1}}
\put(175,38){\framebox(12,12){2}}
\put(175,26) {\framebox(12,12){4}}
\put(187,50){\framebox(12,12){6}}
\put(225,50){\framebox(12,12){3}}
\put(225,38){\framebox(12,12){5}}
\put(237,50){\framebox(12,12){7}}
\put(237,38){\framebox(12,12){8}}
\put(249,50){\framebox(12,12){9}}
\put(160,40){\Huge${(}$}
\put(205,38){,}
\put(270,40){\Huge${)}$}
\put(200,0){Figure 1}
\end{picture} 
\vskip 0in
We say that a box in $(\alpha, \beta)$ has coordinates $(i,j)$ if the box is in the $i$-th row and $j$-th column of either $\alpha$ or $\beta$.  Two boxes in $(\alpha, \beta)$ can have the same coordinates if they occur in the same box in $\alpha$ as in $\beta$; for instance, the left-top-most box in $\alpha$ and the left-top-most box in $\beta$ both have coordinates $(1,1)$. 
\vskip 0in 
The following observation will be used throughout the sequel.
For $m, r_1, n\in \mathbb{N}$, let $m>n$ and $r_1>n$.
Consider a pair of Young diagrams $(\alpha,\beta)$ with $n$ boxes. We construct from the pair $(\alpha,\beta)$ a Young diagram $\mu$ with $mr_1+n$ boxes by adjoining a rectangular box with $r_1$ rows and $m$ columns, as in Figure 2.  The diagram described corresponds to the following partition:
$$\mu = (m+\alpha_1,m+\alpha_2,\ldots,m+\alpha_{r_1},\beta_1,\beta_2,\ldots,\beta_{r_2})$$
\vskip 0in
\noindent
{\bf{Observation:  }} Let $n$, $m$, and  $r_1$ be as above.  Let $[m^{r_1}]\vdash f$.  Then there is a 1-1 correspondence between pairs of Young diagrams with $n$ boxes and Young diagrams containing $[m^{r_1}]$  with $n+f$ boxes.
\vskip 0in
\begin{picture}(400,210)(0,0)
\put(175,200){\line(1,0){100}}
\put(175,100){\line(0,0){100}}
\put(275,100){\line(0,0){100}}
\put(175,100){\line(1,0){100}}
\put(275,200){\line(1,0){75}}
\put(350,175){\line(0,0){25}}
\put(325,175){\line(1,0){25}}
\put(325,150){\line(0,0){25}}
\put(300,150){\line(1,0){25}}
\put(275,125){\line(1,0){25}}
\put(300,125){\line(0,0){25}}
\put(175,25){\line(0,0){75}}
\put(175,25){\line(1,0){25}}
\put(200,25){\line(0,0){25}}
\put(200,50){\line(1,0){25}}
\put(225,50){\line(0,0){25}}
\put(225,75){\line(1,0){25}}
\put(250,75){\line(0,0){25}}
\put(215,210){\vector(-1,0){40}}
\put(235,210){\vector(1,0){40}}
\put(220,210){$m$}
\put(165,160){\vector(0,1){40}}
\put(165,140){\vector(0,-1){40}}
\put(160,150){$r_1$}
\put(300,175){{$\alpha$}}
\put(200,60){{$\beta$}}
\put(180,0){Figure 2}
\put(135,100){$\mu \rightarrow$}
\end{picture}
\vskip 0in
\noindent     
\vskip .1in
\noindent
\subsection{Schur Functions}
\vskip 0in
For details and proofs on the results mentioned in this subsection we refer the reader to \cite{mc} Ch.I, sec. 3.   
Suppose $x_1,x_2,\ldots x_r$ is a finite set of variables. For some partition $\alpha$ we set $x^{\alpha}=x_1^{\alpha_1}\ldots x_r^{\alpha_r}$.  Let $\delta=[r-1,r-2,\ldots,1,0]$.  We define the following determinant. 
$$a_{\alpha+\delta}=\det(x_i^{\lambda_j+r-j})$$ 
This determinant is divisible in $\mathbb{Z}[x_1,\ldots,x_r]$ by the Vandermonde determinant which is given by
$$a_{\delta}=\det(x_i^{r-j}).$$
Then the symmetric \emph{Schur function} is defined as follows: 
$$s_{\alpha}(x_1,\ldots,x_r)={a_{\alpha+\delta}\over{a_{\delta}}}.$$
Take $x_i=q^{i-1}$ ($1\leq i \leq r$).  If $\lambda$ is a partition of length $\leq r$, we have
\begin{equation}
\label{sfq}
s_{\alpha}(1,q,\ldots, q^{r-1})=q^{n(\alpha)}\prod_{1\leq i<j \leq r}{1-q^{\alpha_i-\alpha_j+j-i}\over{1-q^{j-i}}}
\end{equation}
where $n(\alpha)=\sum_{i=1}^{l(\alpha)}{(i-1)\alpha_i}$. We define the following Schur function as a normalization of equation (\ref{sfq}):
\begin{equation}
s_{\alpha,r}(q)={s_{\alpha}(1,q,\ldots,q^{r-1})\over{s_{[1]}(1,q,\ldots,q^{r-1})^{|\alpha|}}}=q^{n(\alpha)}\prod_{1\leq i<j \leq r}{(1-q)(1-q^{\alpha_i-\alpha_j+j-i})\over(1-q^r)(1-q^{j-i})}\notag
\end{equation}
Notice that $s_{\alpha,r}(q)=0$ whenever $l(\alpha)>r$.  Also if $\chi^{\alpha}(q)$ denotes the character of a matrix in $Gl(r)$ with eigenvalues $q^{(-r+1)/2}$, $q^{(-r+3)/2}$, $\ldots$, $q^{(r-1)/2}$ and if the number of boxes in $\alpha$ is $n$, $s_{\alpha,r}(q)=\chi^{\alpha}(q)/(\chi^{[1]}(q))^n$.
\vskip 0in
Let $r=r_1+r_2$; then the Schur function for $\mu=[m+\alpha_1,\ldots,m+\alpha_{r_1},\beta_1,\ldots,\beta_{r_2}]$ is given by 
\begin{equation}
\label{eq:sf}
s_{\mu,r}(q)=q^{n(\mu)}\prod_{1\leq i<j\leq r}{(1-q)(1-q^{\mu_i-\mu_j+j-i})\over{(1-q^r)(1-q^{j-i})}}.
\end{equation}
After some rearrangement of equation~(\ref{eq:sf}) we get the following equation for the Schur function of $\mu$:
\begin{eqnarray}
\label{eq:slr}
s_{\mu,r}(q) &=&  q^{mr_1(r_1-1)/2+r_1|\beta|}{({1-q\over{1-q^r}})^{|\mu|}}s_{\alpha}(1,q,\ldots,q^{r_1-1})s_{\beta}(1,q,\ldots,q^{r_2-1}) \notag \\ 
&&\prod_{i=1}^{r_1}\prod_{j=1}^{r_2}{{1-q^{m+r_1+\alpha_i-\beta_j+j-i}}\over{1-q^{r_1+j-i}}}.
\end{eqnarray} 
The expression of this Schur function in (\ref{eq:slr}) is useful in the proof of the weight formula!
\vskip .5in
\subsection{The Braid Group}
\vskip .1in
\label{sec-prelim}
The braid group of type $A$, $\mathcal{B}_n(A)$, can be defined algebraically by generators $\tilde{\sigma}_1,\tilde{\sigma}_2, \ldots, \tilde{\sigma}_{n-1}$ and relations
\begin{equation}
\tilde{\sigma}_i\tilde{\sigma}_j=\tilde{\sigma}_j\tilde{\sigma}_i \ \ \ \ \mbox{if  $\  |i-j|>1$}, \label{eq:BR1}
\end{equation}
\begin{equation}
\tilde{\sigma}_i\tilde{\sigma}_{i+1}\tilde{\sigma}_i=\tilde{\sigma}_{i+1}\tilde{\sigma}_i\tilde{\sigma}_{i+1} \ \ \ \mbox{if $\  1\leq i\leq n-2$} \label{eq:BR2}
\end{equation}
\indent
Similarly we can define the braid group of type $B$, $\mathcal{B}_n(B)$, by generators $t, \sigma_1,\sigma_2, \ldots, \sigma_{n-1}$ and defining relations given by equations~(\ref{eq:BR1}) and (\ref{eq:BR2}) and
\begin{equation}
\label{eq:BR3}
\sigma_1t\sigma_1t=t\sigma_1t\sigma_1
\end{equation}
\begin{equation}
\label{eq:BR4}
t\sigma_i=t\sigma_i \ \ \ \  \mbox{if  $\ i>1$}. 
\end{equation}
\vskip 0in
Below we illustrate the generator $\tilde{\sigma}_i \in \mathcal{B}_n(A)$, the full-twist $\Delta_3^2 \in \mathcal{B}_3(A)$, and the generators $t, \sigma_i \in \mathcal{B}_n(B)$.
\vskip 0in
\begin{picture}(400,140)(0,0)
\put(40,30){\line(0,0){90}}
\put(38,122){\scriptsize{1}}
\put(45,75){\circle*{1.5}}
\put(50,75){\circle*{1.5}}
\put(55,75){\circle*{1.5}}
\put(65,75){\line(-1,6){7.5}}
\put(75,30){\line(-1,6){6}}
\put(57,122){\scriptsize{i}}
\put(60,30){\line(1,6){15}}
\put(73,122){\scriptsize{i+1}}
\put(80,75){\circle*{1.5}}
\put(85,75){\circle*{1.5}}
\put(90,75){\circle*{1.5}}
\put(95,30){\line(0,0){90}}
\put(95,122){\scriptsize{n}}
\put(60,20){$\tilde{\sigma}_i$}
\put(155,110){\line(-1,5){2}}
\put(151,122){\scriptsize{1}}
\put(165,110){\line(-1,5){2}}
\put(161,122){\scriptsize{2}}
\put(175,109){\line(-1,5){2}}
\put(171,122){\scriptsize{3}}
\put(150,105){\line(6,1){25}}
\put(165,83){\line(-1,5){4}}
\put(175,82){\line(-1,5){4}}
\put(155,83){\line(-1,5){4.3}}
\put(145,78){\line(6,1){30}}
\put(170,60){\line(-1,5){3}}
\put(160,61){\line(-1,5){3}}
\put(150,60){\line(-1,5){3.8}}
\put(145,56){\line(6,1){25}}
\put(150,30){\line(-1,5){5}}
\put(160,30){\line(-1,5){4}}
\put(170,30){\line(-1,5){4}}
\put(155,20){$\Delta_3^2$}
{\linethickness{0.6mm}{\put(230,90){\line(0,0){30}}}}
{\linethickness{0.6mm}{\put(230,30){\line(0,0){50}}}}
\put(240,87){\line(0,0){33}}
\put(237,122){\scriptsize{1}}
\put(216,82){\line(4,1){24}}
\put(216,82){\line(2,-1){10}}
\put(232,75){\line(2,-1){10}}
\put(242,30){\line(0,0){41}}
\put(250,30){\line(0,0){90}}
\put(247,122){\scriptsize{2}}
\put(255,75){\circle*{1.5}}
\put(260,75){\circle*{1.5}}
\put(265,75){\circle*{1.5}}
\put(270,30){\line(0,0){90}}
\put(267,122){\scriptsize{n}}
\put(245,20){$t$}
{\linethickness{0.6mm}{\put(320,30){\line(0,0){90}}}}
\put(330,30){\line(0,0){90}}
\put(327,122){\scriptsize{1}}
\put(335,75){\circle*{1.5}}
\put(340,75){\circle*{1.5}}
\put(345,75){\circle*{1.5}}
\put(355,75){\line(-1,6){7.5}}
\put(340,122){\scriptsize{i}}
\put(350,30){\line(1,6){15}}
\put(353,122){\scriptsize{i+1}}
\put(365,30){\line(-1,6){6}}
\put(370,75){\circle*{1.5}}
\put(375,75){\circle*{1.5}}
\put(380,75){\circle*{1.5}}
\put(385,30){\line(0,0){90}}
\put(383,122){\scriptsize{n}}
\put(350,20){$\sigma_i$}
\put(200,0){Figure 3}
\end{picture}
\vskip 0.05in
In general, the full-twist $\Delta_f^2$ in $f$ strings, is an element in the center of $\mathcal{B}_f(A)$.  Algebraically, $\Delta_f^2 =(\tilde{\sigma}_{f-1}\ldots\tilde{\sigma}_1)^f$. We now define a map from the generators of $\mathcal{B}_n(B)$ into $\mathcal{B}_{f+n}(A)$.  We call the map $\tilde{\rho}_{f,n}$, and we define the image of the generators of $\mathcal{B}_n(B)$ as follows: $\tilde{\rho}_{f,n}(t)=\Delta_f^{-2}\Delta_{f+1}^2$ and $\tilde{\rho}_{f,n}(\sigma_i)=\tilde{\sigma}_{f+i}$ for $i=1,\ldots,n$.  Pictorially, we have the following:
\vskip 0in
\begin{picture}(400,140)(0,0)
\put(120,100){\line(0,0){25}}
\put(130,100){\line(0,0){25}}
\put(135,110){\circle*{1.5}}
\put(140,110){\circle*{1.5}}
\put(145,110){\circle*{1.5}}
\put(150,103){\line(0,0){22}}
\put(162,99){\line(0,0){26}}
\put(108,88){\line(5,1){54}}
\put(120,40){\line(0,0){46}}
\put(130,40){\line(0,0){47}}
\put(108,88){\line(2,-1){10}}
\put(135,60){\circle*{1.5}}
\put(140,60){\circle*{1.5}}
\put(145,60){\circle*{1.5}}
\put(150,40){\line(0,0){50}}
\put(152,69){\line(2,-1){10}}
\put(162,40){\line(0,0){24}}
\put(145,30){\scriptsize{$\Delta^{-2}_f\Delta_{f+1}^2$}}
\put(118,127){\scriptsize{$1$}}
\put(128,127){\scriptsize{$2$}}
\put(147,127){\tiny{f}}
\put(157,127){\tiny{f+1}}
\put(200,40){\line(0,0){85}}
\put(175,40){\line(0,0){85}}
\put(185,60){\circle*{1.5}}
\put(190,60){\circle*{1.5}}
\put(195,60){\circle*{1.5}}
\put(185,110){\circle*{1.5}}
\put(190,110){\circle*{1.5}}
\put(195,110){\circle*{1.5}}
\put(176,127){\tiny{f+2}}
\put(198,127){\tiny{f+n}}
{\linethickness{0.6mm}{\put(12,102){\line(0,0){23}}}}
\put(25,106){\line(0,0){19}}
\put(0,94){\line(2,1){25}}
{\linethickness{0.6mm}{\put(12,40){\line(0,0){55}}}}
\put(0,94){\line(2,-1){10}}
\put(15,87){\line(2,-1){10}}
\put(25,40){\line(0,0){42}}
\put(33,40){\line(0,0){85}}
\put(38,60){\circle*{1.5}}
\put(43,60){\circle*{1.5}}
\put(48,60){\circle*{1.5}}
\put(38,110){\circle*{1.5}}
\put(43,110){\circle*{1.5}}
\put(48,110){\circle*{1.5}}
\put(52,40){\line(0,0){85}}
\put(21,127){\scriptsize{$1$}}
\put(31,127){\scriptsize{$2$}}
\put(50,127){\scriptsize{$n$}}
\put(60,85){\vector(1,0){40}}
\put(75,90){$\tilde{\rho}_{f,n}$}
\put(30,30){\small{$t$}}
{\linethickness{0.6mm}{\put(240,40){\line(0,0){85}}}}
\put(250,40){\line(0,0){85}}
\put(255,85){\circle*{1.5}}
\put(260,85){\circle*{1.5}}
\put(265,85){\circle*{1.5}}
\put(270,124){\line(1,-6){6}}
\put(270,40){\line(1,6){14}}
\put(279,80){\line(1,-6){7}}
\put(284,85){\circle*{1.5}}
\put(289,85){\circle*{1.5}}
\put(294,85){\circle*{1.5}}
\put(299,40){\line(0,0){85}}
\put(310,85){\vector(1,0){30}}
\put(315,90){$\tilde{\rho}_{f,n}$}
\put(248,127){\tiny{1}}
\put(268,127){\tiny{i}}
\put(277,127){\tiny{i+1}}
\put(297,127){\tiny{n}}
\put(270,30){\small{$\sigma_i$}}
\put(350,40){\line(0,0){85}}
\put(355,85){\circle*{1.5}}
\put(360,85){\circle*{1.5}}
\put(365,85){\circle*{1.5}}
\put(365,124){\line(1,-6){6}}
\put(365,40){\line(1,6){14}}
\put(375,80){\line(1,-6){7}}
\put(385,85){\circle*{1.5}}
\put(390,85){\circle*{1.5}}
\put(395,85){\circle*{1.5}}
\put(404,40){\line(0,0){85}}
\put(402,127){\tiny{f+n}}
\put(348,127){\tiny{1}}
\put(373,127){\tiny{f+i+1}}
\put(358,127){\tiny{f+i}}
\put(375,30){\small{$\tilde{\sigma}_{f+i}$}}
\put(200,5){Figure 4}
\end{picture}
\vskip 0in
\begin{proposition}
\label{prop:rep1}
Let $n,f\in \mathbb{N}$ and $\tilde{\rho}_{f,n}$ be as defined above for the generators of $\mathcal{B}_n(B)$. Then $\tilde{\rho}_{f,n}$ is a well-defined group homomorphism.
\end{proposition}
\proof
It suffices to prove that $\tilde{\rho}_{f,n}$ preserves the relations of $\mathcal{B}_n(B)$. It is convenient to express $\tilde{\rho}_{f,n}(t)=\tilde{\sigma}_f\ldots\tilde{\sigma}_2\tilde{\sigma}_1^2\tilde{\sigma}_2\ldots\tilde{\sigma}_f$.  Relations (\ref{eq:BR1}),(\ref{eq:BR2}) and (\ref{eq:BR4}) follow immediately from the definition of $\tilde{\rho}_{f,n}$ and the definition of $\mathcal{B}_n(A)$.  To show that (\ref{eq:BR3}) holds we will use the commuting property of the full-twist.  
$$
\sigma_{f+1}\Delta_{f+1}^2\Delta_f^{-2}\sigma_{f+1}\Delta_{f+1}^2\Delta_f^{-2} = \Delta_{f+2}^2\Delta_{f+1}^{-2}\Delta_{f+1}^2\Delta_f^{-2} \\ = \Delta_{f+1}^{2}\Delta_f^{-2}\sigma_{f+1}\Delta_{f+1}^{2}\Delta_f^{-2}\sigma_{f+1}.
$$
It might be easier to see that this relation holds from the picture, see figure 5.
This completes the proof.~$\Box$ 
\vskip .08in
\begin{picture}(400,150)(0,0)
\put(80,117){\line(0,0){28}}
\put(80,87){\line(0,0){21}}
\put(80,50){\line(0,0){28}}
\put(85,140){\circle*{1.5}}
\put(90,140){\circle*{1.5}}
\put(95,140){\circle*{1,5}}
\put(100,121){\line(0,0){24}}
\put(100,90){\line(0,0){22}}
\put(100,50){\line(0,0){31}}
\put(110,123){\line(0,0){22}}
\put(120,121){\line(0,0){24}}
\put(65,110){\line(5,1){55}}
\put(85,100){\circle*{1.5}}
\put(90,100){\circle*{1.5}}
\put(95,100){\circle*{1.5}}
\put(110,88){\line(0,0){26}}
\put(65,110){\line(2,-1){10}}
\put(112,97){\line(2,-1){10}}
\put(122,50){\line(0,0){42}}
\put(65,80){\line(5,1){45}}
\put(85,65){\circle*{1.5}}
\put(90,65){\circle*{1.5}}
\put(95,65){\circle*{1.5}}
\put(65,80){\line(2,-1){10}}
\put(102,67){\line(2,-1){10}}
\put(112,50){\line(0,0){12}}
\put(137,50){\line(0,0){95}}
\put(142,90){\circle*{1.5}}
\put(147,90){\circle*{1.5}}
\put(152,90){\circle*{1.5}}
\put(157,50){\line(0,0){95}}
\put(78,147){\tiny{1}}
\put(98,147){\tiny{f}}
\put(104,147){\tiny{f+1}}
\put(118,147){\tiny{f+2}}
\put(135,147){\tiny{f+3}}
\put(155,147){\tiny{f+n}}
\put(57,40){\scriptsize{$\tilde{\sigma}_{f+1}\Delta_f^{-2}\Delta_{f+1}^2\tilde{\sigma}_{f+1}\Delta^{-2}_f\Delta_{f+1}^2$}}
\put(180,90){\Huge{$=$}}
\put(230,117){\line(0,0){28}}
\put(235,140){\circle*{1.5}}
\put(240,140){\circle*{1.5}}
\put(245,140){\circle*{1,5}}
\put(250,121){\line(0,0){24}}
\put(260,119){\line(0,0){26}}
\put(270,86){\line(0,0){59}}
\put(215,110){\line(5,1){45}}
\put(230,83){\line(0,0){24}}
\put(235,100){\circle*{1.5}}
\put(240,100){\circle*{1.5}}
\put(245,100){\circle*{1.5}}
\put(250,87){\line(0,0){24}}
\put(215,110){\line(2,-1){10}}
\put(252,97){\line(2,-1){10}}
\put(274,50){\line(0,0){8}}
\put(215,75){\line(5,1){55}}
\put(230,50){\line(0,0){23}}
\put(235,65){\circle*{1.5}}
\put(240,65){\circle*{1.5}}
\put(245,65){\circle*{1.5}}
\put(250,50){\line(0,0){26}}
\put(215,75){\line(2,-1){10}}
\put(264,63){\line(2,-1){10}}
\put(262,87){\line(0,0){4}}
\put(262,50){\line(0,0){29}}
\put(287,50){\line(0,0){95}}
\put(292,90){\circle*{1.5}}
\put(297,90){\circle*{1.5}}
\put(302,90){\circle*{1.5}}
\put(307,50){\line(0,0){95}}
\put(228,147){\tiny{1}}
\put(248,147){\tiny{f}}
\put(254,147){\tiny{f+1}}
\put(268,147){\tiny{f+2}}
\put(283,147){\tiny{f+3}}
\put(306,147){\tiny{f+n}}
\put(215,40){\scriptsize{$\Delta_f^{-2}\Delta_{f+1}^2\tilde{\sigma}_{f+1}\Delta^{-2}_f\Delta_{f+1}^2\tilde{\sigma}_{f+1}$}}
\put(175,5){Figure 5}
\end{picture}
\begin{cor}
The representations $\rho_{f,n}$ can be extended in a natural way to representations of the corresponding braid group algebras.
\end{cor}

\section{Hecke Algebra of type $B_n$}
\indent
In this context, we will mean by the Hecke algebra $H_n(q,Q)$ of type $B_n$ the free complex algebra with 1 and generators $t, g_1, g_2, \ldots, g_{n-1}$  and parameters $q,Q \in \mathbb{C}$ with defining relations
\begin{enumerate}
\item[(H1)] \ \ \ $g_{i}g_{i+1}g_{i}=g_{i+1}g_{i}g_{i+1}$ \ \ \ 
                        for $i=1,2, \ldots, n-2$;
\item[(H2)] \ \ \  $g_{i}g_{j}= g_{j}g_{i}$,  \ \ \ 
                        whenever $\mid i-j \mid \geq 2$;
\item[(H3)] \ \ \   $g_i^2 = (q-1)g_i +q$ \ \ \ \ for $i=1,2, \ldots, n-1$;
\item[(H4)] \ \ \   $t^2 = (Q-1)t +Q$;
\item[(H5)] \ \ \   $tg_1tg_1=g_1tg_1t$;
\item[(H6)] \ \ \   $tg_i=g_it$ \ \ \ \ \ for $i \geq 2$.
\end{enumerate}
Note that for $q=1$ and $Q=1$, (H3) and (H4) become $g_i^2=1$ and $t^2=1$, respectively.  In this case the generators satisfy exactly the same relations as a set of simple reflections of the hyperoctahedral group, $\mathcal{H}_n$.  It is known that $H_n(q,Q)\cong \mathbb{C}\mathcal{H}_n$ (the group algebra of $\mathcal{H}_n$)  if $q$ is  not a root of unity and $Q\neq -q^s$ for $-n<s<n$. $\mathbb{C}\mathcal{H}_n$ is semisimple since it is a finite complex group algebra. This implies that $H_n(q,Q)$ is semisimple for generic values of $q$ and $Q$.  Furthermore the simple modules are in 1-1 correspondence with double partitions. Their decomposition rule and their dimensions are the same as for $\mathcal{H}_n$. [See Bourbarki, Groups et algebres de Lie IV, V, VI]
\vskip 0in
\indent
Hoefsmit, in \cite{h}, has written down explicit irreducible representations of $H_n(q,Q)$ for each pair of Young diagrams.  This demonstrates that dimension over $\mathbb{C}(q,Q)$ of $H_n(q,Q)$ is $2^nn!$. 
 Let $V_{(\alpha,\beta)}$  be an $H_n(q,Q)$ module then we have the following restriction rule:
\begin{equation}
V_{(\alpha,\beta)}{\Big|}_{H_{n-1}(q,Q)}\cong \bigoplus_{(\alpha,\beta)'\subset(\alpha,\beta)}V_{(\alpha,\beta)'},
\end{equation}
where $(\alpha,\beta)'$ is the double partition obtained from $(\alpha,\beta)$ by subtracting one from one of the parts. This isomorphism follows from the bijection $t^{(\alpha,\beta)}\rightarrow t^{(\alpha,\beta)'}$ of $T_{(\alpha,\beta)}$ and $\bigcup_{(\alpha,\beta)'\subset (\alpha,\beta)}T_{(\alpha,\beta)'}$, where $t^{(\alpha,\beta)'}$ is the pair of standard tableaux obtained by removing the box containing $n$. Moreover, $\pi^{(n)}=\bigoplus_{(\alpha,\beta)\vdash n}\pi_{(\alpha,\beta)}$ is a faithful representation of $H_n(q,Q)$.
\vskip 0in
\indent
Note that (H1) and (H2) are the defining relations for the braid group $\mathcal{B}_n(A)$.  So representations of $H_n(q,Q)$ also yield in particular representations of $\mathcal{B}_n(A)$. Obviously we also obtain representations of the braid group of type $B$.
\vskip 0in
\indent
The Hecke algebras satisfies the following embedding of algebras  $H_0\subset H_1 \subset H_2 \subset \cdots $ We shall denote 
$$ H_{\infty}(q,Q):= \bigcup_{n \geq 0}{H_n(q,Q)}$$
\vskip 0.02in
Observe that (H3) and (H4) imply that $t$ and $g_i$ have at most 2 eigenvalues each, hence also at most 2 projections corresponding to these eigenvalues.
\vskip 0.1in
\noindent
\emph{Right and Double Coset Representatives}
\vskip 0.05in  
The fact that $q$ is invertible in $\mathbb{C}(q,Q)$ implies that the generators $g_i$ are also invertible in $H_n(q,Q)$.  In fact, the inverse of the generators is given by
$$g_i^{-1} = q^{-1}g_i + (q^{-1}-1)1 \in H_n(q,Q)$$
This implies that the following element is well-defined in $H_n(q,Q)$: $t_i^{'}=g_i\cdots g_1tg_1^{-1} \cdots g_i^{-1}$. 
We use this elements to define the set $\mathcal{D}_n$ as a subset of $H_n(q,Q)$.  If $n=1$, we let $\mathcal{D}_1 := \{1,t\}$.  For $n\geq 2$ we have
$$\mathcal{D}_n:=\{1,g_{n-1},t^{'}_{n-1}\}.$$
Also the set of right coset representatives of $H_{n-1}(q,Q)$ in $H_n(q,Q)$ is given as follows:
$$\mathcal{R}_n :=\{ 1,\ t'_{n-1},\  g_{n-1}g_{n-2}\ldots g_{n-k},\  g_{n-1}g_{n-2}\ldots g_{n-k}t'_{n-k-1} \ (1\leq k\leq n-1)\}.$$
Note that $\mathcal{R}_1 = \{1,t\}$.  The elements $\{ r_1r_2\ldots r_n| r_i \in \mathcal{R}_i\}$ form a $\mathbb{C}(q,Q)$-basis of $H_n(q,Q)$, see \cite{gl}.

\section{Representations of the Hecke algebra of type $B$ onto a Reduced Hecke Algebra of type $A$}
\vskip 0in
In this section we will show that for a specialization of the Hecke algebra of type $B$, we have a homomorphism onto a reduced Hecke algebra of type $A$. Before we prove this assertion we introduce some necessary background.
\vskip 0in
In this section we assume that for a fixed integer $n$ and $m,r_1\in \mathbb{N}$ we have $m>n$ and $r_1>n$. For these integers we choose $\lambda=[m^{r_1}]$.  Then as we observed in Section 1, there is a one-to-one correspondence between double partitions of $n$ and partitions of $n+mr_1$ containing $\lambda$.  Once we choose $\lambda$ we fix a standard tableau $t^{\lambda}$.
We also assume throughout that $q$ is not a root of unity.  
\vskip 0in
The representations of $H_n(q,Q)$ defined by Hoefsmit, (see \cite{h}), depend on rational functions with denominators $(Qq^d+1)$ where $d\in \{0,\pm 1,\ldots, \pm (n-1)\}$.  So if $Q=-q^{r_1+m}$ then $1-q^{r_1+m+d}\neq 0$ as long as $d\neq -(r_1+m)$, which implies that all representations of $H_n(q,-q^{r_1+m})$ are well-defined as long as $r_1+m>n$.  Thus the specialized algebra, $H_n(q,-q^{r_1+m})$ is well-defined and semisimple.
\vskip 0in
\subsection{The Hecke Algebra of type $A$}
By the Hecke algebra of type $A_{n-1}$, $H_n(q)$, we mean the free complex algebra with generators $\tilde{g}_1$, \ldots, $\tilde{g}_{n-1}$ and 1 with parameter $q\in \mathbb{C}$ and  defining relations
\begin{enumerate}
\item[$(\widetilde{H1})$] \ \ \ $\tilde{g}_{i}\tilde{g}_{i+1}\tilde{g}_{i}=\tilde{g}_{i+1}\tilde{g}_{i}\tilde{g}_{i+1}$ \ \ \ 
                        for $i=1,2, \ldots, n-2$;
\item[$(\widetilde{H2})$] \ \ \  $\tilde{g}_{i}\tilde{g}_{j}= \tilde{g}_{j}\tilde{g}_{i}$,  \ \ \ 
                        whenever $\mid i-j \mid \geq 2$;
\item[$(\widetilde{H3})$] \ \ \   $\tilde{g}_i^2 = (q-1)\tilde{g}_i +q$ \ \ \ \ for $i=1,2, \ldots, n-1$;
\end{enumerate} 
$H_n(q)$ is semisimple whenever $q$ is not a root of unity.  If $\mu \vdash n$, then $(\pi_{\mu}, V_{\mu})$ denotes the irreducible representation of $H_n(q)$ indexed by $\mu$. Here $V_{\mu}$ is the vector space with orthonormal basis given by $\{v_{t^{\mu}}\}$ where $t^{\mu}$ is a standard tableau of shape $\mu$. These representations can be considered as $q$-analogs of Young's orthogonal representations of the Symmetric group (see \cite{h} or \cite{w1}). We now describe these representations in more detail.
\vskip 0in  
Let $d(t^{\mu},i)=c(i+1)-c(i)-r(i+1)+r(i)$ where $c(i)$ and $r(i)$ denote respectively, the column and row of the box containing the number $i$, and for $d\in \mathbb{Z}\setminus\{0\}$ let 
$$a_d(q)={q^d(1-q)\over{1-q^d}} \ \ \   \mbox{and} \ \ \  c_d(q)={\sqrt{(1-q^{d+1})(1-q^{d-1})}\over{1-q^d}}.$$
Then one defines $\pi_{\mu}$ on the vector space $V_{\mu}$  by 
$$\pi_{\mu}(\tilde{g}_{i})v_{t^{\mu}}=a_d(q)v_{t^{\mu}}+c_d(q)v_{\tilde{g}_i(t^{\mu})}$$
where $d=d(t^{\mu},i)$ and $\tilde{g}_i(t^{\mu})$ is the tableau obtained by interchanging $i$ and $i+1$ in $t^{\mu}$.  This representation follows from the ones defined in \cite{w1} by setting $\pi_{\mu}(\tilde{g}_i)=(q+1)\pi_{\mu}(e_i)-Id_{V_{\mu}}$, where $e_i$ is the eigenprojection corresponding to the characteristic value $q$ of $\tilde{g}_i$.  Moreover $\pi_n=\bigoplus_{\mu\vdash n}\pi_{\mu}$ is a faithful representation of $H_n(q)$.  Thus we have that $H_n(q)\cong \bigoplus_{\mu\vdash n}\pi_{\mu}(H_n(q))$, where $\pi_{\mu}(H_n(q))\cong M_{f^{\mu}}$ and $ M_{f^{\mu}}$ is the ring of complex $f^{\mu}\times f^{\mu}$ matrices, and $f^{\mu}$ is the number of standard tableaux of shape $\mu$.
\vskip 0in
If $\mu$ is a partition of $n$ then denote by $\mu'$ the partition obtained by decreasing one part of $\mu$ by one; and $\mu''$ is obtained by decreasing one part of $\mu'$ by one.
Let $t^{\mu}$ be a standard  tableau of shape $\mu$, then $t^{\mu}\rightarrow (t^{\mu})'$ defines a bijection between $T_{\mu}$ and $\bigcup_{\mu'\subset \mu}T_{\mu'}$, where $(t^{\mu})'$ is the Young tableau obtained by removing the box containing $n$ from $t^{\mu}$.  This in turn yields the following isomorphism of modules
\begin{equation}
\label{eq:modA}
V_{\mu}{\Big|}_{H_{n-1}(q)}\cong\bigoplus_{\mu'\subset\mu}V_{\mu'}.
\end{equation} 
\vskip 0in
A special element in $H_f(q)$ is the \emph{full-twist} defined algebraically by $\Delta_f^2:=(\tilde{g}_{f-1}\ldots \tilde{g}_1)^f$.  The full-twist is an element in the center of $H_f(q)$.  The action of $\Delta_f^2$ on $(\pi_{\nu}, V_{\nu})$ is described in the following lemma, the proof of which appears in \cite{w2}, pg. 261.  
\begin{lemma}
\label{prop:ftr}
Let $\nu \vdash f$.
Then the full-twist acts by a scalar $\alpha_{\nu}$ on the irreducible representation $(\pi_{\nu},V_{\nu})$ of $H_f(q)$ where 
$$\alpha_{\nu}=q^{f(f-1)-\sum_{i<j}(\nu_i+1)\nu_j}.$$
\end{lemma}
\subsection{A Homomorphism of $H_n(q,-q^{r_1+m})$ onto a reduced algebra of $H_{n+f}(q)$}
In (\cite{w1}, Cor. 2.3) Wenzl defined a special set of minimal idempotents of $H_f(q)$ indexed by the standard tableaux. The sum of these idempotents is 1.  These minimal idempotents are well-defined whenever $H_f(q)$ is semisimple. 
\vskip 0in
If $p\in H_f(q)$ is an idempotent then the \emph{reduced algebra} of $H_f(q)$ is $pH_f(q)p:= \{pap\, |\, a\in H_f(q)\}$.
Let $\lambda\vdash f$ and $t^{\lambda}$ be a standard tableau of shape $\lambda$. Then $p_{t^{\lambda}}$ is  the minimal idempotent indexed by $t^{\lambda}$. Thus, the reduced algebra decomposes as follows:
$$p_{t^{\lambda}}H_{n+f}(q)p_{t^{\lambda}} \cong \bigoplus_{\lambda \subset \mu,\atop \mu\vdash n+f}\pi_{\mu}(p_{t^{\lambda}})\pi_{\mu}(H_{n+f}(q))\pi_{\mu}(p_{t^{\lambda}}).$$ 
Notice that if $\mu$ is a partition of $n+f$ and $\mu$ does not contain $\lambda$ then $\pi_{\mu}(p_{t^{\lambda}})$ is the zero matrix.  In particular, if we 
choose  $\lambda$ to be rectangular, we have that the reduced algebra $p_{t^{\lambda}}H_{f+1}(q)p_{t^{\lambda}}$ has only two nonzero irreducible modules indexed by partitions $[m+1, m^{r_1-1}]$ and $[m^{r_1},1]$, since these are the only partitions of $f+1$ which contain $\lambda$. 
\vskip 0in
Recall that in Section 1 we showed that there is a homomorphism $\tilde{\rho}_{f,n}$ from the braid group $\mathcal{B}_n(B)$ into the braid group $\mathcal{B}_{f+n}(A)$.  In what follows we will show that it is possible to extend $\tilde{\rho}_{f,n}$ to a homomorphism of the corresponding Hecke algebras.
\begin{lemma}
\label{lemma:hom}
For a fixed integer $n$ and $m,r_1\in \mathbb{N}$, let $m>n$ and $r_1>n$. Set $f=mr_1$ and assume $\lambda=[m^{r_1}]$ and $\gamma=[m^{r_1},1]$. Let $\alpha_{\lambda}$ and $\alpha_{\gamma}$ be as in Lemma ~\ref{prop:ftr}.  Choose a minimal idempotent $p_{t^{\lambda}}$ in $H_f(q)$. We define a map $\rho_{f,n}$ for the generators of $H_n(q,-q^{r_1+m})$ as follows:  
$$\rho_{f,n}(1)=p_{t^{\lambda}},\hspace{20pt}\rho_{f,n}(t)=-{\alpha_{\lambda}\over{\alpha_{\gamma}}}p_{t^{\lambda}}\Delta_{f}^{-2}\Delta_{f+1}^2\hspace{20pt} \mbox{and} \hspace{20pt} \rho_{f,n}(g_i)=p_{t^{\lambda}}\tilde{g}_{i+f}$$
for $i=1,\ldots, n-1$, (see Figure 6 for a pictorial definition). Then $\rho_{f,n}$ extends to a well-defined homomorphism of algebras, $\rho_{f,n}:H_n(q,-q^{r_1+m})\rightarrow p_{t^{\lambda}}H_{n+f}(q)p_{t^{\lambda}}$. 
\end{lemma}
\proof It is enough to check that $\rho_{f,n}$ preserves the relations of the Hecke algebra of type $B$.  First notice that since $p_{t^{\lambda}}\in H_f(q)$, it commutes with all $\tilde{g}_j$ for all $j>f$. From this observation and the defining relations of $H_{n+f}(q)$ it follows that we only need to check the relations involving $\rho_{f,n}(t)$. For relations (H5) and (H6), (see pg. 7), notice that $p_{t^{\lambda}}$ commutes with the full-twist  $\Delta_{f+1}^2$, $\Delta_f^{-2}$ and with $\tilde{g}_{f+1}$. Thus these two relations will hold by the proof of Prop. 1.2.
\vskip 0in
It remains to be shown that $\rho_{f,n}(t)$ has two eigenvalues: $-1$ and $Q=-q^{r_1+m}$. In particular, we want to show that $\rho_{f,n}(t)$ acts by a scalar on two irreducible modules of $H_{f+1}(q)$ and by $0$ on all others.  Since $\lambda$ is rectangular this means that there are only two partitions of $f+1$ containing it. 
\vskip 0in
By lemma \ref{prop:ftr} we have that the full-twist, $\Delta_{f+1}^2$, acts by the scalar $\alpha_{\beta}$ (resp. $\alpha_{\gamma}$) on the irreducible module indexed by $\beta$ (resp. $\gamma$).  Also $\Delta_f^{-2}$ acts by $\alpha_{\lambda}^{-1}$ on the irreducible modules of $H_{f+1}(q)$ which are indexed by Young diagrams which contain $\lambda$. Therefore, we have that $\Delta_f^{-2}\Delta_{f+1}^2$ acts by $\alpha_{\beta}\alpha_{\lambda}^{-1}$ on the module $V_{\beta}$ and by $\alpha_{\gamma}\alpha_{\lambda}^{-1}$ on the module $V_{\gamma}$.  Therefore, we have that $-{\alpha_{\lambda}\over{\alpha_{\gamma}}}p_{t^{\lambda}}\Delta_f^{-2}\Delta_{f+1}^2$ acts by $-1$ on $V_{\gamma}$ and by $-{\alpha_{\beta}\over{\alpha_{\gamma}}}$ on $V_{\beta}$, and by zero on all other modules, since $p_{t^{\lambda}}$ kills all modules which do not contain $\lambda$.  
\vskip 0in
In order to determine the constant $-{\alpha_{\beta}\over{\alpha_{\gamma}}}$  we only need to substitute the partitions  $\beta=[m+1, m^{r_1-1}]$ and $\gamma=[m^{r_1},1]$ in the formula given in lemma \ref{prop:ftr}. Thus we obtain
$$-{\alpha_{\beta}\over{\alpha_{\gamma}}}=-q^{-\sum_{i<j}(\beta_i+1)\beta_j+\sum_{i<j}(\gamma_i+1)\gamma_j}=-q^{r_1+m}.$$
This concludes this proof. $\Box$
\vskip 0.03in
The following figure gives a pictorial definition of the homomorphism $\rho_{f,n}$ defined in the previous lemma. 
\vskip 0in
\vskip 0in
\begin{picture}(400,130)(0,0)
{\linethickness{0.6mm}{\put(12,92){\line(0,0){23}}}}
\put(25,96){\line(0,0){19}}
\put(0,84){\line(2,1){25}}
{\linethickness{0.6mm}{\put(12,30){\line(0,0){55}}}}
\put(0,84){\line(2,-1){10}}
\put(15,77){\line(2,-1){10}}
\put(25,30){\line(0,0){42}}
\put(33,30){\line(0,0){85}}
\put(38,70){\circle*{1.5}}
\put(43,70){\circle*{1.5}}
\put(48,70){\circle*{1.5}}
\put(52,30){\line(0,0){85}}
\put(21,117){\scriptsize{$1$}}
\put(31,117){\scriptsize{$2$}}
\put(50,117){\scriptsize{$n$}}
\put(60,75){\vector(1,0){40}}
\put(75,80){$\rho_{f,n}$}
\put(119,105){\framebox(23,10){$p_{t^{\lambda}}$}}
\put(120,85){\line(0,0){20}}
\put(120,40){\line(0,0){36}}
\put(118,117){\scriptsize{$1$}}
\put(125,70){\circle*{1.5}}
\put(130,70){\circle*{1.5}}
\put(135,70){\circle*{1.5}}
\put(139,90){\line(0,0){15}}
\put(139,40){\line(0,0){40}}
\put(138,117){\tiny{f}}
\put(119,30){\framebox(23,10){$p_{t^{\lambda}}$}}
\put(152,87){\line(0,0){28}}
\put(108,78){\line(5,1){44}}
\put(108,78){\line(2,-1){10}}
\put(142,59){\line(2,-1){10}}
\put(152,30){\line(0,0){24}}
\put(146,117){\tiny{f+1}}
\put(164,30){\line(0,0){85}}
\put(161,117){\tiny{f+2}}
\put(168,70){\circle*{1.5}}
\put(173,70){\circle*{1.5}}
\put(178,70){\circle*{1.5}}
\put(183,30){\line(0,0){85}}
\put(181,117){\tiny{f+n}}
{\linethickness{0.6mm}{\put(220,30){\line(0,0){85}}}}
\put(227,30){\line(0,0){85}}
\put(225,117){\tiny{1}}
\put(232,75){\circle*{1.5}}
\put(237,75){\circle*{1.5}}
\put(242,75){\circle*{1.5}}
\put(242,114){\line(1,-6){6}}
\put(251,70){\line(1,-6){7}}
\put(240,117){\tiny{i}}
\put(242,30){\line(1,6){14}}
\put(249,117){\tiny{i+1}}
\put(255,75){\circle*{1.5}}
\put(260,75){\circle*{1.5}}
\put(265,75){\circle*{1.5}}
\put(274,30){\line(0,0){85}}
\put(272,117){\tiny{n}}
\put(280,75){\vector(1,0){30}}
\put(285,80){$\rho_{f,n}$}
\put(319,105){\framebox(20,10){$p_{t^{\lambda}}$}}
\put(320,40){\line(0,0){65}}
\put(319,117){\tiny{1}}
\put(324,75){\circle*{1.5}}
\put(329,75){\circle*{1.5}}
\put(334,75){\circle*{1.5}}
\put(338,40){\line(0,0){65}}
\put(336,117){\tiny{f}}
\put(319,30){\framebox(20,10){$p_{t^{\lambda}}$}}
\put(348,30){\line(0,0){85}}
\put(344,117){\tiny{f+1}}
\put(353,75){\circle*{1.5}}
\put(358,75){\circle*{1.5}}
\put(363,75){\circle*{1.5}}
\put(365,114){\line(1,-6){6}}
\put(375,70){\line(1,-6){7}}
\put(358,117){\tiny{f+i}}
\put(365,30){\line(1,6){14}}
\put(373,117){\tiny{f+i+1}}
\put(383,75){\circle*{1.5}}
\put(388,75){\circle*{1.5}}
\put(393,75){\circle*{1.5}}
\put(398,30){\line(0,0){85}}
\put(397,117){\tiny{f+n}}
\put(200,5){Figure 6}
\end{picture}
\vskip 0in
We have shown that $\rho_{f,n}$ is a homomorphism. In the remainder of this section we will show that it is onto.  In particular, we will show that the irreducible representations of the reduced algebra are also irreducible representations of $H_n(q,-q^{r_1+m})$.
\vskip 0in
Notice that for $\mu \vdash n+f$ there is a 1-1 correspondence between standard skew tableaux of shape $\mu/\lambda$ and tableaux $t^{\mu}$ which contain $t^{\lambda}$. For this reason we denote by $T_{\mu/\lambda}$ the set of standard tableaux which contain $t^{\lambda}$. Notice that the order of $T_{\mu/\lambda}$ is equal to the number of standard skew tableaux of shape $\mu/\lambda$.  For our choice of $m$, i.e. $m>n$, we have that $\mu/\lambda$ will consist of two parts which can be interpreted as a double partition, say $(\delta,\gamma)$.  In this case, the order of $T_{\mu/\lambda}$  is ${n\choose |\delta|}f^{\delta}f^{\gamma}$, where $f^{\gamma}$ is the number of standard tableaux of shape $\gamma$; this formula is given by Hoefsmit in \cite{h}.       
\vskip 0in
Define $V_{\mu/\lambda}$ as the complex vector space with orthonormal basis $\{v_{t^{\mu}}\,|\, t^{\mu}\in T_{\mu/\lambda}\}$. Notice that $V_{\mu/\lambda}$ is a subspace of  $V_{\mu}$ and $p_{t^{\lambda}}V_{\mu}=V_{\mu/\lambda}$, where $V_{\mu}$ has basis indexed by standard tableaux of shape $\mu$. 
\vskip 0in
\noindent
{\bf{Observation:  }} Let $\mu\vdash n+f$ and $V_{\mu}$ be an irreducible $H_{n+f}(q)$ module.  Then $p_{t^{\lambda}}V_{\mu}$ is an irreducible $p_{t^{\lambda}}H_{n+f}(q)p_{t^{\lambda}}$ module.
\vskip 0in
This observation implies that there is a set of irreducible representations of $p_{t^{\lambda}}H_{n+f}(q)p_{t^{\lambda}}$ indexed by Young diagrams with $n+f$ boxes, which contain $\lambda$. It follows from equation~(\ref{eq:modA}) that 
\begin{equation}
\label{eq:sum}
p_{t^{\lambda}}V_{\mu}{\Big|}_{p_{t^{\lambda}}H_{n+f-1}(q)p_{t^{\lambda}}}\cong\bigoplus_{\lambda\subset\mu'\subset \mu}V_{\mu'/\lambda}
\end{equation}
where $\mu'$ has $n+f-1$ boxes and $p_{t^{\lambda}}V_{\mu'}=V_{\mu'/\lambda}$.
\vskip 0in
We let $z_{\mu}$ denote the minimal central idempotents of $H_{n+f}(q)$, and $z_{\mu/\lambda}$ denote the minimal central idempotents of $p_{t^{\lambda}}H_{n+f}(q)p_{t^{\lambda}}$. Notice that $z_{\mu/\lambda}$ has rank equal to the number of skew standard tableaux of shape $\mu/\lambda$.
\begin{theorem}
\label{th:repba}
 Let $f, n$ be as in Lemma \ref{lemma:hom} and assume that q is not a root of unity. Then $\rho_{f,n}$ as defined in Lemma \ref{lemma:hom} is an onto homomorphism. 
\end{theorem}
\proof  In Lemma \ref{lemma:hom} we showed that $\rho_{f,n}$ is a homomorphism.  Thus it only remains to show that it is onto.  The proof is by induction on $n$.  For $n=1$, we have $\rho_{f,1}:H_1(q,-q^{r_1+m})\rightarrow p_{t^{\lambda}}H_{f+1}(q)p_{t^{\lambda}}$. Since $\lambda\vdash f$ is a rectangular diagram there are only two Young diagrams with $f+1$ boxes which contain $\lambda$, i.e. $[m+1, m^{r_1-1}]$ and $[m^{r_1}, 1]$.  As we showed in Lemma \ref{lemma:hom} the action of $\rho_{f,1}(t)$ on the representation indexed by $[m+1, m^{r_1-1}]$ (resp. $[m^{r_1}, 1]$) is $-q^{r_1+m}$ (resp. $-1$).  And both representations are 1 dimensional.      
\vskip 0in
The algebra $H_1(q,-q^{r_1+m})$ has two irreducible representations indexed by $([1], \emptyset)$ and $(\emptyset, [1])$. Both of these representations are 1 dimensional and $t\in H_1(q,-q^{r_1+m})$ acts by a scalar on these representations. The action of $t$ on $V_{([1], \emptyset)}$ (resp. $V_{(\emptyset,[1])}$) is $-q^{r_1+m}$ (resp. $-1$). Since $q$ is not a root of unity, these representations are irreducible and nonequivalent. 
This shows that $\pi_{([1], \emptyset)}\cong \pi_{[m+1,m^{r_1-1}]}$ and $\pi_{(\emptyset,[1])}\cong \pi_{[m^{r_1},1]}$. 
\vskip 0in
Assume that for $n>1$ we have $\rho_{f,n}$ is onto. Then for all $\mu$ containing $\lambda$, $V_{\mu/\lambda}$ is an irreducible $H_{n+1}(q,-q^{r_1+m})$, and $V_{\mu/\lambda}{\Big|}_{H_n(q,-q^{r_1+m})}\cong \bigoplus_{\lambda\subset\mu'\subset \mu}V_{\mu'/\lambda}$, as in equation ~(\ref{eq:sum}). 
\vskip 0in
The remainder of this proof is similar to the proof of irreducibility of modules of the Hecke algebra of type $A$ in \cite{w1}, Theorem 2.2.  
\vskip 0in
Let $\mu\vdash n+f+1$ be as described above.
By the induction assumption, $H_n(q, -q^{r_1+m})$ is a semisimple algebra with minimal central idempotents $z_{\mu'/\lambda}$.  Let $0\neq W\subset V_{\mu/\lambda}$ be an  $H_{n+1}(q, -q^{r_1+m})$  module.  But $V_{\mu/\lambda}$ decomposes as an $H_n(q, -q^{r_1+m})$ module into the direct sum of irreducible modules $V_{\mu'/\lambda}$ since $\rho_{f,n}$ is onto $p_{t^{\lambda}}H_{n+f}(q)p_{t^{\lambda}}$.  Thus, there exists a $\mu'\vdash$ $n+f$ such that $V_{\mu'/\lambda}\subset W$.  Let $\tilde{\mu}'\neq \mu'$ be another Young diagram with $n+f$ boxes such that $\lambda\subset \tilde{\mu}'\subset \mu$. There is exactly one $\mu''\vdash$ $n+f-1$ contained in both $\mu'$ and $\tilde{\mu}'$ such that $\mu''$ contains $\lambda$.  Let $t^{\mu}\in T_{\mu/\lambda}$ be such that $(t^{\mu})' \in T_{\mu'/\lambda}$ and $(t^{\mu})''\in T_{\mu''/\lambda}$.  Then $(\tilde{g}_{n+f}(t^{\mu}))'\in T_{\tilde{\mu}'/\lambda}$ and therefore 
$$\pi_{\mu}(z_{\tilde{\mu}'/\lambda}) \pi_{\mu}(\tilde{g}_{n+f})v_{t^{\mu}}=c_dv_{\tilde{g}_{n+f}(t^{\mu})}\in V_{\tilde{\mu}'/\lambda}$$
where $d$ is the axial distance in $t^{\mu}$ between $n+f$ and $n+f+1$.  Since $q$ is not a root of unity then $c_d$ is well-defined and nonzero.  Hence the irreducible $H_n(q, -q^{r_1+m})$ module, $V_{\tilde{\mu}'/\lambda}$, is contained in W.  But $\tilde{\mu}'$ was arbitrary, therefore $W\supset\bigoplus_{\mu'\subset \mu} V_{\mu'}$.
\vskip 0in
Next we show that the $V_{\mu/\lambda}$ are mutually nonisomorphic $H_{n+1}(q,-q^{r_1+m})$ modules.  As we observed above this is true for $n=1$. For $n=2$ there are five irreducible modules;  4 are one dimensional and 1 is two dimensional.  We must check that the one dimensional modules are nonequivalent.  By the definition of the action of $t$ and $\tilde{g}_{f+1}$ we have that $t$ acts by $-q^{r_1+m}$ on $V_{[m+2,m^{r_1-1}]}$ and $V_{[m+1,m+1,m^{r_1-2}]}$.  But $\tilde{g}_{f+1}$ acts by $q$ on $V_{[m+2,m^{r_1-1}]}$ and by $-1$ on $V_{[m+1,m+1,m^{r_1-2}]}$. In a similar way we can show that $V_{[m^{r_1},2]}$ and $V_{[m^{r_1},1^2]}$ are nonequivalent.  And since $t$ acts by $-1$ on these last two modules, we have that they are nonequivalent to the former two.   If $n>2$ and $\mu$ and $\tilde{\mu}$ are two distinct partitions of $n+f+1$ which contain $\lambda$, then there exists a $\mu'\supset\lambda$ such that $\mu'\subset\mu$ but $\mu'\not\subset{\tilde{\mu}}$. The proof of this fact is found in \cite{w1} Lemma 2.11.  We would also like to remark that restricting to Young diagrams containing $\lambda$ does not affect the result.  Hence, $V_{\mu/\lambda}$ and $V_{\tilde{\mu}/\lambda}$ differ already as $H_n(q,-q^{r_1+m})$ modules.~$\Box$
\vskip 0.03in
\remark  Notice that this theorem formalizes the idea we indicated in Section 1 about associating pairs of Young diagrams $(\alpha, \beta)$ with one Young diagram where we adjoin the $m\times r_1$ rectangle, see Figure 2. 
\vskip 0.03in
It is well-known that there exits a duality between the quantum group $U_q(sl(r))$  and the Hecke algebra of type $A$.  This duality is the quantum analogue of the  Schur-Weyl duality between the general linear group, $GL(n)$, and the symmetric group, $S_n$, (see \cite{d}, \cite{ji1}, and \cite{ji2}).
\vskip 0in
The following is an easy Corollary of Theorem \ref{th:repba}. 
\begin{cor}
The action of the specialized Hecke algebra of type $B$, $H_n(q,-q^{r_1+m})$ and the diagonal action of $U_q(sl(r))$ on $V_{\lambda}\otimes V^{\otimes n}$ have the double centralizing property in $End(V_{\lambda}\otimes V^{\otimes n})$.
\end{cor}
The proof of this corollary follows immeditely from the duality between $U_q(sl(r))$ and the Hecke algebra of type $A$.
\vskip .03in
\remark  To relate the above to the literature we make the following remark.  However, this remark will not be used in the sequel. For definitions we refer the reader to~\cite{ji1}.
\vskip 0in 
The results of this section imply that there is an R-matrix representation of $H_n(q,-q^{r_1+m})$ on $V_{\lambda}\otimes V^{\otimes n}$, where $V$ is the fundamental module of $U_q(sl(r))$ and $V_{\lambda}$ is the irreducible module corresponding to $\lambda$.  In this R-matrix representation $t$ acts on $V_{\lambda}\otimes V=V_{\beta}\oplus V_{\gamma}$ ($\gamma$ and $\beta$ as defined above) as a scalar.  And $g_i\rightarrow R_i$ is given by the same R-matrices as for the Hecke algebra of type $A$ where $R_i$ acts on the $i$-th and $i+1$ copies of $V$.
\vskip 0in
\noindent
 \section{\bf {Markov Traces on the Hecke Algebra of type $B$}}
In this section we give the necessary background on Markov traces.  We refer the reader to \cite{gl} or \cite{g} for the details.
\vskip 0in
A \emph{trace} function on $H_{\infty}(q,Q)$ is a $\mathbb{C}(q,Q)$-linear map $\phi:H_{\infty}(q,Q)\longrightarrow \mathbb{C}(q,Q)$ such that $\phi (hh')=\phi (h'h)$ for all $h,h' \in H_{\infty}(q,Q)$.
This definition is in fact valid for any associative algebra over a commutative ground ring.  In the case of the group algebra it is clear that every trace function is constant on the conjugacy classes of the underlying group. Notice that $\phi(hh')-\phi(h'h)=0$; this means that $\phi(hh'-h'h)=\phi([h,h'])=0$, which implies that the commutators are in the kernel of a trace function in an algebra. 
\vskip 0in
The weights we are going to give in this paper correspond to a trace that satisfies the following definition.
\vskip 0.05in
\noindent
{\bf{Definition:  }}
 Let $z\in \mathbb{C}(q,Q)$ and $tr : H_{\infty}(q,Q) \longrightarrow \mathbb{C}(q,Q)$ be an $\mathbb{C}(q,Q)$-linear map.  Then $tr$ is called a \emph{Markov trace} (with parameter $z$) if the following conditions are satisfied:
\vskip 0in
(1) $tr$ is a trace function on $H_{\infty}(q,Q)$;
\vskip 0in
(2) $tr(1)=1$ (normalization);
\vskip 0in
(3) $tr(hg_n)= ztr(h)$ for all $n \geq 1$ and $h \in H_n(q,Q)$.
\vskip .05in
The name of these traces comes from their invariance under the Markov moves for closed braids. Remember that the Hecke Algebra is a quotient of the braid group. We note that all generators $g_i$ (for $i = 1,2, \ldots$) are conjugate in $H_{\infty}(q,Q)$.  In particular, any trace function on $H_{\infty}(q,Q)$ must have the same value on these elements.  This explains why the parameter $z$ is independent of $n$ in rule (3) of this definition.
\vskip 0in
Geck and Pfeiffer in \cite{gp} showed that a trace function on the Hecke algebra is uniquely determined by its value on basis elements corresponding to representatives of minimal length in the various conjugacy classes of the underlying Coxeter group.  Also  representatives of minimal length in the classes of Coxeter groups of classical types are of the form $d_1 \ldots d_n$, where $d_i$ is a distinguished double coset representative of the $H_i(q,Q)$ with respect to $H_{i-1}(q,Q)$. 
\vskip 0in
Let $tr$ be a Markov trace with parameter $z$, and let $d_i\in \mathcal{D}_i$ for $i=1,\ldots,n$.  Then 
$$tr(d_1\ldots d_n)=z^a tr(t'_0t'_1\ldots t'_{b-1})$$
where $a$ is the number of factors $d'_i$ which are equal to $g_{i-1}$ and $b$ is the number of factors which are equal to $t'_{i-1}$.  Thus, $tr$ is uniquely determined by its parameter $z$ and the values on the elements in the set $\{t'_0t'_1\ldots t'_{i-1}\,|\,i=1,2,\ldots\}$.
\vskip 0in 
Conversely, given $z,\, y_1, \, y_2, \ldots \in \mathbb{C}(q,Q)$ then there exist a unique Markov trace on $H_{\infty}(q,Q)$ such that 
$tr(t'_0t'_1\ldots t'_{k-1})=y_k$ for all $k\geq 1$.  For details on these results see~\cite{gl}, Theorem 4.3.
\vskip 0in
We are particularly interested in the special case when  $y_i=y^i$ for all $i\in \mathbb{N}$, this is the case when there is only two parameters.  In this case, we have that if  $d_i$ is a distinguished double coset representative then $tr(d_ix)=\xi tr(x)$ where $\xi=y$ or $z$. The proof of the following proposition is found in \cite{gl}.
\begin{proposition}
Let $z, y \in \mathbb{C}(q,Q)$ and $tr: H_{\infty}(q,Q) \longrightarrow \mathbb{C}(q,Q)$ be a Markov trace with parameter $z$ such that $tr(t^{'}_0t^{'}_1 \ldots t^{'}_{k-1}) =y^{k}$ for all $k \geq 1$ then
\vskip 0in
\noindent 
\centerline{$tr(ht^{'}_{n,0})= ytr(h)$ for all $n \geq 0$ and $h\in H_n(q,Q)$}
\vskip 0in
\noindent
where $t'_{n,0}=g_n\ldots g_1tg_1^{-1}\ldots g_n^{-1}$ or $g_n^{-1}\ldots g_1^{-1}tg_1\ldots g_n$.
\end{proposition}     
\indent
We point out that the converse is also true. 
A proof of the existence of such traces for every $z\in \mathbb{C}$ is given in \cite{gl}.  We will compute the weight vectors for this Markov trace on $H_n(q,Q)$.

\section{\bf{The Weight Formula}}

In this section we define for every pair of partitions, $(\alpha, \beta)$, a rational function in $q$ and $Q$, $W_{(\alpha,\beta)}(q,Q)$. We will show that this function is a weight for the Markov trace defined by Geck and Lambropolou \cite{gl} for the Hecke algebra of type $B$.  If we denote the weights by $\omega_{(\alpha,\beta)}$ then the Markov trace, $tr$, can be written as follows:
\begin{equation}
tr(x)=\sum_{(\alpha,\beta)\vdash n}{\omega_{(\alpha,\beta)}\chi^{(\alpha,\beta)}(x)},
\end{equation}
where $ x \in H_n(q,Q)$ and $\chi^{(\alpha,\beta)}$ is the character (the usual trace) of the irreducible representation of $H_n(q,Q)$ indexed by $(\alpha,\beta)$.
\vskip 0in
First we define a rational function in $q$ and $Q$ for any arbitrary double partition $(\alpha,\beta)$ such that $l(\alpha)\leq r_1$ and $l(\beta)\leq r_2$. If $l(\alpha)=s< r_1$ then $\alpha_i=0$ for $i=s+1,\ldots r_1$, similarly for $\beta$. Let $r=r_1+r_2$. 
\begin{eqnarray}
\label{eq:wd}
W_{(\alpha,\beta)}(q,Q) & = & 
q^{n(\alpha)+n(\beta)}({{1-q}\over{1-q^r}})^{|\alpha|+|\beta|}\prod_{1\leq i<j\leq r_1}{{1-q^{\alpha_i-\alpha_j+j-i}} 
\over{1-q^{j-i}}}
\prod_{1\leq i<j\leq r_2} {{1-q^{\beta_i-\beta_j+j-i}} \over{1-q^{j-i}}} \notag
\\
&&\times  \prod_{i=1}^{r_1}\prod_{j=1}^{r_2}{{Qq^{\alpha_i-i}+q^{\beta_j-j}}\over
{Qq^{-i}+q^{-j}}} 
\end{eqnarray}
Notice that this function can be expressed as a product of Schur functions 
\begin{equation}
\label{eq:wsf}
W_{(\alpha,\beta)}(q,Q)= q^{r_1|\beta|}{s_{\alpha}(1,q,\ldots,q^{r_1-1})
s_{\beta}(1,q,\ldots,q^{r_2-1})\over{s_{[1]}(1,q,\ldots, q^{r_1+r_2-1})^{|\alpha|+|\beta|}}}\prod_{i=1}^{r_1}
\prod_{j=1}^{r_2}{(1+Qq^{\alpha_i-\beta_j+j-i})\over(1+Qq^{j-i})}.
\end{equation}
Recall that  $s_{\alpha}(1,q,\ldots,q^{r_1-1})=\prod_{1\leq i<j\leq r_1}{{1-q^{\alpha_i-\alpha_j+j-i}} \over{1-q^{j-i}}}$ is the symmetric Schur function defined in Section 1, see equation~(\ref{sfq}).
\vskip 0in
\noindent
{\bf{Observation: }} Let $1-r_1\leq s \leq r_2-1$. Assume that $q$ is not a root of unity and $Q\neq -q^{-s}$.  Then $W_{(\alpha,\beta)}(q,Q)$ is an analytic rational function.
\vskip 0in
The rectangular Young diagram with $r$ rows which has $m$ boxes in the first $r_1$ rows and $0$ boxes in the remaining rows, i.e. $[m^{r_1}]$, has Schur function given by 
$$s_{[m^{r_1}],r}(q)={q^{mr_1(r_1-1)/2}\prod_{i=1}^{r_1}\prod_{j=1}^{r_2}{1-q^{m+r_1+j-i}\over{1-q^{r_1+j-i}}}\over{s_{[1]}(1,q,\ldots,q^{r-1})^{mr_1}}}.$$
\vskip .1in
We are going to assume that for fixed integer $n$ and $m, r_1 \in \mathbb{N}$, we have $m>n$ and $r_1>n$. Then for any double partition of $n$, $(\alpha,\beta)$, we have $\mu=[m+\alpha_1,\ldots,m+\alpha_{r_1},\beta_1,\ldots,\beta_{r_2}]$ is a partition of $n+mr_1$.  Then by equation~(\ref{eq:slr}) from Section 1 we have the following equation
\begin{equation}
{s_{\mu,r}(q)\over{s_{[m^{r_1}],r}(q)}}=q^{r_1|\beta|}{s_{\alpha}(1,q,\ldots,q^{r_1-1})s_{\beta}(1,q,\ldots q^{r_2-1})\over{s_{[1]}(1,q,\ldots, q^r)^{|\alpha|+|\beta|}}}\prod_{i=1}^{r_1}\prod_{j=1}^{r_2}{1-q^{m+r_1+\alpha_i-\beta_j+j-i}\over{1-q^{m+r_1+j-i}}}. \notag
\end{equation}
Observe that we have the following equality: 
\begin{equation}
\label{eq:lim}
{s_{\mu,r}(q)\over{s_{[m^{r_1}],r}(q)}}=W_{(\alpha,\beta)}(q,-q^{r_1+m}).
\end{equation}
Notice that $W_{(\alpha,\beta)}(q,-q^{r_1+m})$ is well-defined since $r_1+m>r_1$, and as we observed before $W_{(\alpha,\beta)}(q,Q)$ is undefined for $Q=-q^{-s}$ where $1-r_1<s<r_2-1$.  So $W_{(\alpha,\beta)}(q,-q^{r_1+m})$ is an analytic rational function.
\begin{lemma}
$$
W_{(\alpha,\beta)}(q,Q)=\sum_{(\alpha,\beta)\subset(\gamma,\eta)} W_{(\gamma,\eta)}(q,Q)
$$
where $(\alpha,\beta)\subset(\gamma,\eta)$ means that $(\gamma,\eta)$ is obtained by adding one box to $(\alpha,\beta)$.
\end{lemma}
\proof  Assume $\mu=[m+\alpha_1,\ldots,m+\alpha_{r_1},\beta_1,\ldots,\beta_{r_2}]$. 
By the Littlewood-Richardson rule  for Schur functions (see \cite{mc}) we have the following: 
$$s_{[1],r}s_{\mu,r}=\sum_{\mu\subset\nu, \atop |\nu|=|\mu|+1}s_{\nu,r}$$
where $s_{[1],r}(q)=1$.
Now divide both sides of this equation by $s_{[m^{r_1}],r}(q)$ and we get by equation~(\ref{eq:lim}).
\begin{equation}
\label{eq:*}
W_{(\alpha,\beta)}(q,-q^{r_1+m})=\sum_{(\alpha,\beta)\subset(\gamma,\eta)}W_{(\gamma,\eta)}(q,-q^{r_1+m}).
\end{equation}
Since $W_{(\alpha,\beta)}(q,Q)$ is an analytic rational function and the above equation holds for all $Q=-q^{r_1+m}$, thus we have that 
\begin{equation}
W_{(\alpha,\beta)}(q,Q)=\sum_{(\alpha,\beta)\subset (\gamma,\eta)}W_{(\gamma,\eta)}(q,Q).
\end{equation}
holds for all values of $Q$.$\Box$
\vskip .1in
In \cite{w1} Wenzl showed that the weights of the Markov trace (with parameter $z=q^r{(1-q)\over{(1-q^r)}}$) on the Hecke algebra of type $A$  are given by the symmetric Schur function  $s_{\mu,r}(q)$ as described in Section 1. 
\vskip 0in
Let $\lambda=[m^{r_1}]\vdash f$.  Since we assumed that $m>n$ and $r_1>n$, then we have that for all $\mu\vdash n+f$ we have that $\mu/\lambda$ can be interpreted as a double partition $(\alpha, \beta)$ of $n$. 
Now we fix $t^{\lambda}$ a standard tableau of shape $\lambda$.  Then the reduced algebra $p_{t^{\lambda}}H_{n+f}(q)p_{t^{\lambda}}$ has $p_{t^{\lambda}}$ as the identity.  The Markov trace for the reduced algebra is given by the renormalized Markov trace of $H_{n+f}(q)$.  By renormalization we mean that we must divide the trace by the trace of the identity, i.e. $tr(p_{t^{\lambda}})=s_{\lambda,r}(q)$.  
Therefore, we have that ${s_{\mu,r}(q)\over{s_{\lambda,r}(q)}}$ are the weights of $p_{t^{\lambda}}H_{n+f}(q)p_{t^{\lambda}}$.  Notice that this implies that $W_{(\alpha,\beta)}(q,Q)$ reduce to the weights for the reduced algebra when $Q=-q^{r_1+m}$. 
\begin{lemma} Let $g_{n-1}\in H_n(q,Q)$, $z={q^r(1-q)\over{(1-q^r)}}$ and $W_{(\alpha,\beta)}(q,Q)$ as defined in equation (\ref{eq:wd}). Then for any $x\in H_n(q,Q)$
\begin{equation}
\label{eq:def}
tr(x)=\sum_{(\alpha,\beta)\vdash n}W_{(\alpha,\beta)}(q,Q)\chi^{(\alpha,\beta)}(x)
\end{equation}
defines a well-defined trace which satisfies the Markov property, i.e. $tr(hg_{n-1})=ztr(h)$, where $h\in H_{n-1}(q,Q)$.
\end{lemma} 
\proof It is clear that $tr$ is indeed a trace.  We must show that it satisfies the Markov property.  We have assumed that $\mu/\lambda=(\alpha,\beta)$, then $W_{(\alpha,\beta)}(q,-q^{r_1+m})={s_{\mu,r}(q)\over{s_{\lambda,r}(q)}}$. We also have that $\chi^{\mu/\lambda}=\chi^{(\alpha,\beta)}\Big |_{Q=-q^{r_1+m}}$ since we have shown in section 3 that $V_{\mu/\lambda}$ is an irreducible module of $H_n(q,-q^{r_1+m})$.  Then we have 
$$\sum_{\mu\vdash n+f}{s_{\mu,r}(q)\over{s_{\lambda,r}(q)}}\chi^{\mu/\lambda}(x)=\sum_{(\alpha,\beta)\vdash n}W_{(\alpha,\beta)}(q,-q^{r_1+m})\chi^{(\alpha,\beta)}(x)$$
defines a Markov trace for the reduced algebra $p_{t^{\lambda}}H_{n+f}(q)p_{t^{\lambda}}$ with parameter $z={q^r(1-q)\over{(1-q^r)}}$. But in Theorem \ref{th:repba} we showed that $\rho_{f,n}:H_n(q,-q^{r_1+m})\rightarrow p_{t^{\lambda}}H_{n+f}(q)p_{t^{\lambda}}$ is an onto homomorphism. In particular we showed that the irreducible modules of $p_{t^{\lambda}}H_{n+f}(q)p_{t^{\lambda}}$ are irreducible modules for $H_n(q,-q^{r_1+m})$.  Therefore, we have that the above equation also defines a trace which satisfies the Markov property for $H_n(q,-q^{r_1+m})$. We know that $W_{(\alpha,\beta)}(q,Q)$ and $\chi^{(\alpha,\beta)}$ are analytic functions; thus by the identity theorem in complex analysis, since the weights work for all $Q=-q^{r_1+m}$,  they must work for all $Q$.~$\Box$ 
\vskip .1in
Lemma 5.1 and Lemma 5.2 imply that the function $W_{(\alpha,\beta)}(q,Q)$ defined in equation (\ref{eq:wd}) is a weight function for a Markov trace with parameter $z=q^r(1-q)/(1-q^r)$.  In Section 4 we noted that a Markov trace on the Hecke algebra of type $B$ is uniquely determined by a parameter $z$ and by the values on the set $\{t_0't_1'\ldots t_{k-1}'\, |\,k\geq 1\}$.  Therefore, we must still compute the values of the $tr$ on this set in order to completely determine this trace.  To compute these values, we need the following definitions.  
\vskip 0.1in
Let $A\subset B$ be a pair of semisimple finite algebras.  Given a trace, $tr$, nondegenerate on both $A$ and $B$, the \emph{conditional expectation} $E_A: B\rightarrow A$ is defined by setting $E_A(b)$ to be the unique element of $A$ such that $$tr(E_A(b)a)=tr(ba)\ \mbox{for all} \  a\in A \mbox{ and } b\in B$$
$E_A$ is well-defined and unique.
\vskip 0in
Let $B$ be represented via left multiplication on itself, where we write $L^2(B, tr)$ to denote the representation space $B$ to distinguish it from the algebra $B$.  We use $b_{\xi}$ to denote the elements in $L^2(B, tr)$.  
\vskip 0in
One obtains from $E_A$ an idempotent $e_A: L^2(B,tr)\rightarrow L^2(B,tr)$ defined by $e_Ab_{\xi}=E_A(b)_{\xi}$.  The idempotent $e_A$ can be thought of as an orthogonal projection onto $A$ with respect to the bilinear form $(b_{\xi},c_{\xi})\rightarrow tr(bc)$. The algebra $\langle B,e_A\rangle$ generated by $B$, acting via left multiplication on $L^2(B,tr)$, and by the idempotent $e_A$ is called the Jones \emph{basic construction} for $A\subset B$.  For the proof of the following results see \cite{j1}.
\begin{theorem}
Let $A$, $B$, $E_A$, $e_A$, and $tr$ be as defined above. Then
\vskip 0in
(a) The algebra $\langle B,e_A\rangle$ is isomorphic to the centralizer $End_A B$ of $A$ acting by left multiplication on $B$.  In particular, it is semisimple.
\vskip .02in
(b) There is a 1-1 correspondence between the simple components of $A$ and $End_A B$ such that if $p\in A_i$ is a minimal idempotent, $pe_A$ is a minimal idempotent of $\langle B,e_A \rangle_i$.  
\vskip .02in 
(c) $e_Abe_A=E_A(b)e_A$ for all $b\in B$.
\end{theorem}   
\vskip 0in
In our case we have the pair of semisimple algebras $H_{n+f-1}(q)\subset H_{n+f}(q)$. The corresponding orthogonal projection is $p_{[1^r]}$, where $r$ is the length of the partitions indexing the irreducible representations.
It might also be helpful to keep the following figure in mind, since it clearly shows the veracity of the next lemma.
\vskip 0in
\begin{picture}(400,160)(0,0)
\put(24,110){\framebox(27,15){$ p_{t^{\lambda}}$}}
\put(25,90){\line(0,0){20}}
\put(30,90){\line(0,0){20}}
\put(50,95){\line(0,0){15}}
\put(14,84){\line(6,1){84}}
\put(14,84){\line(2,-1){10}}
\put(25,37){\line(0,0){8}}
\put(30,37){\line(0,0){8}}
\put(50,37){\line(0,0){8}}
\put(50,60){\line(0,0){23}}
\put(50,125){\line(0,0){20}}
\put(48,147){\tiny{f}}
\put(24,45){\framebox(27,15){$ p_{t^{\lambda}}$}}
\put(60,98){\line(0,0){47}}
\put(60,76){\line(0,0){12}}
\put(60,37){\line(0,0){30}}
\put(57,147){\tiny{f+1}}
\put(67,80){\circle*{1.5}}
\put(72,80){\circle*{1.5}}
\put(77,80){\circle*{1.5}}
\put(85,98){\line(0,0){47}}
\put(85,72){\line(0,0){17}}
\put(85,37){\line(0,0){27}}
\put(75,147){\tiny{f+n-1}}
\put(55,74){\line(4,-1){43}}
\put(25,60){\line(0,0){21}}
\put(30,60){\line(0,0){21}}
\put(35,80){\circle*{1.5}}
\put(40,80){\circle*{1.5}}
\put(45,80){\circle*{1.5}}
\put(95,125){\framebox(37,14){$p_{[1^r]}$}}
\put(98,98){\line(0,0){27}}
\put(98,50){\line(0,0){14}}
\put(98,139){\line(0,0){6}}
\put(95,147){\tiny{f+n}}
\put(110,50){\line(0,0){75}}
\put(110,139){\line(0,0){6}}
\put(115,80){\circle*{1.5}}
\put(120,80){\circle*{1.5}}
\put(125,80){\circle*{1.5}}
\put(130,50){\line(0,0){75}}
\put(130,139){\line(0,0){6}}
\put(123,147){\tiny{f+n+r-1}}
\put(175,75){\Huge{=}}
\put(95,35){\framebox(37,14){$p_{[1^r]}$}}
\put(24,23){\framebox(62,14){$h$}}
\put(25,125){\line(0,0){20}}
\put(24,147){\tiny{1}}
\put(30,125){\line(0,0){20}}
\put(29,147){\tiny{2}}
\put(243,125){\framebox(40,14){$p_{[1^r]}$}}
\put(244,50){\line(0,0){75}}
\put(244,139){\line(0,0){6}}
\put(242,147){\tiny{1}}
\put(251,50){\line(0,0){75}}
\put(251,139){\line(0,0){6}}
\put(250,147){\tiny{2}}
\put(255,80){\circle*{1.5}}
\put(260,80){\circle*{1.5}}
\put(265,80){\circle*{1.5}}
\put(270,50){\line(0,0){75}}
\put(270,139){\line(0,0){6}}
\put(267,147){\tiny{r-1}}
\put(280,139){\line(0,0){6}}
\put(280,95){\line(0,0){30}}
\put(280,50){\line(0,0){20}}
\put(279,147){\tiny{r}}
\put(243,35){\framebox(40,14){$p_{[1^r]}$}}
\put(295,110){\framebox(25,12){$ p_{t^{\lambda}}$}}
\put(298,75){\line(0,0){35}}
\put(298,57){\line(0,0){13}}
\put(298,122){\line(0,0){23}}
\put(293,147){\tiny{r+1}}
\put(303,90){\circle*{1.5}}
\put(308,90){\circle*{1.5}}
\put(313,90){\circle*{1.5}}
\put(318,78){\line(0,0){32}}
\put(318,57){\line(0,0){15}}
\put(318,122){\line(0,0){23}}
\put(314,147){\tiny{r+f}}
\put(280,95){\line(3,-1){15}}
\put(320,83){\line(2,-1){10}}
\put(280,70){\line(6,1){50}}
\put(295,45){\framebox(25,12){$ p_{t^{\lambda}}$}}
\put(335,36){\line(0,0){109}}
\put(330,147){\tiny{r+f+1}}
\put(340,90){\circle*{1.5}}
\put(345,90){\circle*{1.5}}
\put(350,90){\circle*{1.5}}
\put(298,35){\line(0,0){9}}
\put(318,35){\line(0,0){9}}
\put(297,23){\framebox(62,13){$h$}}
\put(355,36){\line(0,0){109}}
\put(354,147){\tiny{r+f+n-1}}
\put(160,0){Figure 7}
\end{picture}
\vskip 0in
Before stating the next lemma, we would like to define a tensor product of Hecke algebras given in \cite{gw} by Goodman and Wenzl.  We have $H_n(q)\otimes H_m(q)\subset H_{n+m}(q)$ defined by $a\otimes b=a(shift_n(b))$ for $a\in H_n(q)$ and $b\in H_m(q)$, where $shift_n$ is the operator which sends $g_i$ to $g_{i+n}$ for all $i$.  In particular, this tensor product allows us to multiply minimal idempotents using the generalized Littlewood-Richardson rule in \cite{gw}.  Denote by $1_n$ the identity in $H_n(q)$.
\begin{lemma}
A Markov trace on the Hecke Algebra of type $A$ induces a Markov trace on the Hecke algebra of type B, which satisfies the condition that $tr(t'_nx)=ytr(x)$, where $x\in H_n(q,Q)$.  In particular, this implies that $y_k=y^k$ for all $k\in\mathbb{N}$. 
\end{lemma}
\proof Let $tr$ be a Markov trace for $H_{n+f}(q)$ and $tr_{p_{t^{\lambda}}}$ be the Markov trace corresponding to the reduced algebra.  We have shown that the weights for the reduced algebra define a Markov trace for $H_n(q,-q^{r_1+m})$.  
\vskip 0in
Assume that $E_{H_{n+f-1}(q)}:H_{n+f}(q)\rightarrow H_{n+f-1}(q)$ is the unique conditional expectation with respect to the Markov trace defined by the weight function. Thus, $E_{H_{n+f-1}(q)}(\tilde{g}_{n+f}h)=zh$ for all $h\in H_{n+f-1}(q)$ and $z=q^r{(1-q)\over{(1-q^r)}}$.
\vskip 0in
We denote the image of $t_n'$ under the homomorphism $\rho_{n+f}$ by $p_{t^{\lambda}}\tau_{n+f}$.
Consider the element $p_{t^{\lambda}}\tau_{n+f}h$ where $h\in H_{n+f-1}(q)$.  We would like to compute the conditional expectation of 
$p_{t^{\lambda}}\tau_{n+f}h$.  Observe that $$E_{H_{n+f-1}(q)}(p_{t^{\lambda}}\tau_{n+f}h)=E_{H_{n+f-1}(q)}(p_{t^{\lambda}}\tau_{n+f})h.$$  
Therefore it suffices to compute the conditional expectation for $p_{t^{\lambda}}\tau_{n+f}$.  Consider the following 
$$(p_{t^{\lambda}}\otimes 1_{n-1} \otimes p_{[1^r]})(\tau_{n+f} \otimes 1_{r-1})(p_{t^{\lambda}}\otimes 1_{n-1}\otimes p_{[1^r]})(h\otimes 1_r)$$
The above expresion is equal to the left hand side of the following equation (see figure 7). 
$$[((p_{[1^r]}\otimes p_{t^{\lambda}})(1_{r-1}\otimes \tau_{n+f})(p_{[1^r]}\otimes p_{t^{\lambda}}))\otimes 1_{n-1}](1_r \otimes h) = (const)(p_{[1^r]}\otimes p_{t^{\lambda}}\otimes 1_{n-1})(1_r \otimes h)$$
Since we are restricted to $r$ rows we have by the generalization of the Littlewood-Richardson rule, see \cite{gw}, that $p_{[1^r]}\otimes p_{t^{\lambda}}$ is a minimal idempotent.
But Theorem 5.3 (c) implies that $(const)p_{t^{\lambda}}=E_{H_{n+f-1}}(q)(\tau_{n+f}p_{t^{\lambda}})$. Thus,
$$tr_{p_{t^{\lambda}}}(\tau_{n+f}h)=(const) tr_{p_{t^{\lambda}}}(h).$$
This implies that  we have a Markov trace with the property $tr(t_n'h)=(const)tr(h)$ for all $h\in H_{n+f-1}(q)$. But this trace is also a trace for $H_n(q,-q^{r_1+m})$.  Since the above property holds for all $Q=-q^{r_1+m}$, then it must hold for all $Q$.  Thus our assertion is proved.~$\Box$
\begin{theorem}
If $tr$ is a Markov trace on the Hecke algebra of type $B$, with parameter $z=q^r(1-q)/(1-q^{r_1+r_2})$, such that $tr(t_0't_1'\ldots, t'_{k-1})=y^k$ for $k\geq1$.  Then the weights are given by  $W_{(\alpha,\beta)}(q,Q)$ as defined in equation~(\ref{eq:wd}) with $y={(q^{r_2}Q+1)(1-q^{r_1})/(1-q^r}) -1$.
\end{theorem}
\proof  The fact that the weights define a Markov trace follows from Lemmas 5.1,  5.2 and 5.4.  It remains to show that the weight formula is indeed given by this values of $y={(Qq^{r_2}+1)(1-q^{r_1})/(1-q^{r_1+r_2})} -1$. By Lemma 5.1 it suffices to compute $tr(t)$ in $H_1(q)$. This is a straight forward computation,
we have $W_{([1],\emptyset)}(q,Q)={(1-q^{r_1})(1+Qq^{r_2})\over{(1-q^r)(1+Q)}}$ and $W_{(\emptyset,[1])}(q,Q)={q^{r_1}(1-q^{r_1})(1+Qq^{-r_1})\over{(1-q^r)(1+Q)}}$. Also $\chi^{([1],\emptyset)}(t)=Q$ and $\chi^{(\emptyset,[1])}(t)=-1$.  Using these values we compute 
$$tr(t)=Q W_{([1],\emptyset)}(q,Q)-W_{(\emptyset,[1])}(q,Q)={(Qq^{r_2}+1)(1-q^{r_1})\over{(1-q^r)}}-1. \hskip 70pt \Box$$

\section{\bf{Markov Trace for the Hecke algebra of type D}}
\vskip 0in
The easiest way to study Markov traces on the Hecke algebras of type $D$ is by embedding these algebras into those for type $B$, and then applying the results for type $B$. 
\vskip 0in
As observed by Hoefsmit \cite{h}, in order to obtain an embedding of $H_n^D(q)$ into $H_n^B(q,Q)$ we have to set the parameter $Q$ equal to 1.  Notice that in this case we have $t^2=1$.  The Hecke algebra of type $D$ is generated by $u=tg_1t,\ g_1, \ldots, g_{n-1}$ satisfying the following relations:
\begin{enumerate}
\item[(D1)] $g_ig_{i+1}g_i=g_{i+1}g_ig_{i+1}$ for $i=1,\ldots n-2$;
\item[(D2)] $g_ig_j=g_jg_i$ whenever $|i-j|\geq 2$;
\item[(D3)] $g_i^2=(q-1)g_i+q$ for all $i$;
\item[(D4)] $g_iu=ug_i$ for all $i$;
\item[(D5)] $ u^2=(q-1)u+q$.
\end{enumerate}
We have $H_n^D(q)\subset H_n^B(q,Q)$ for all $n$; then $H_{\infty}^D(q)=\bigcup_{n>1}H_n^D(q)\subset H_{\infty}^B(q,Q)$.  Geck in \cite{g} shows that the restriction of a Markov trace on $H_{\infty}^B(q,Q)$ is a Markov trace on $H_{\infty}^D(q)$ and both have the same parameter.  Furthermore, he shows that every Markov trace on $H_{\infty}^D(q)$ can be obtained in this way. 
\vskip 0in
From Hoefsmit \cite{h} we know that the simple components for $H_n^D(q)$ are indexed by double partitions $(\alpha,\beta)$.  If $\alpha\neq \beta$  we have that the $H_n^B(q,Q)$-modules $V_{(\alpha,\beta)}$ and $V_{(\beta,\alpha)}$ are simple, equivalent $H_n^D(q)$-modules.  And if $\alpha = \beta$ we have that the $H_n^B(q,Q)$-module $V_{(\alpha,\alpha)}$ decomposes into two simple nonequivalent $H_n^D(q)$-modules, i.e. $V_{(\alpha,\alpha)_i}$ with $i=1,2$. Using Bratteli diagrams we have the following relations for simple modules of the $H_n^B(q,Q)$ and $H_n^D(q)$
\vskip 0in
\begin{picture}(400,120)
\put(10,70){type $B$}
\put(60,100){$\alpha \neq \beta$}
\put(60,70){$(\alpha,\beta)$}
\put(140,70){$(\beta,\alpha)$}
\put(75,65){\line(2,-1){40}}
\put(115,45){\line(2,1){40}}
\put(100,35){$(\alpha,\beta)$}
\put(10,35){type $D$}
\put(250,100){$\alpha=\beta$}
\put(275,70){$(\alpha,\alpha)$}
\put(240,35){$(\alpha,\alpha)_1$}
\put(255,45){\line(2,1){40}}
\put(294,65){\line(2,-1){40}}
\put(320,35){$(\alpha,\alpha)_2$}
\put(180,0){figure 8}
\end{picture}
\vskip 0in  
\begin{proposition}
The weight formula for the Markov trace on the Hecke algebra of type $D$ with parameters $z=q^r{(1-q)\over{(1-q^r)}}$ and $y={(Qq^{r_2}+1)(1-q^{r_1})\over{(1-q^r)}}-1$ is given as follows: 
$$W_{(\alpha,\beta)}^D(q)=W_{(\alpha,\beta)}(q,1)+W_{(\alpha,\beta)}(q,1),\ \ \mbox{ if } \alpha \neq \beta$$
and
$$W_{(\alpha,\alpha)_i}^D(q)=W_{(\alpha,\alpha)}(q,1), \mbox{ for } i=1,2 \ \  \mbox{ if }\  \alpha = \beta$$
where $W_{(\alpha,\beta)}(q,1)$ denote the weight evaluated at $Q=1$ of the Hecke algebra of type $B$.  
\end{proposition}
\proof  The proof of this theorem follows directly from the inclusion matrix of the Hecke algebra of type $D$ into the Hecke algebra of type $B$. Recall that in order to obtain the weight vector for the Hecke algebra of type $D$ we multiply the inclusion matrix for $H_n^D(q)\subset H_n^B(q,Q)$ with the weight vector for type $B$. $\Box$ 
\vskip .1in
\bibliographystyle{plain}

University of California, San Diego, 
\vskip 0in \noindent
Department of Mathematics, La Jolla, CA 93106, USA
\vskip 0in \noindent
\emph{E-mail address:} rorellan@@math.ucsd.edu 
\end{document}